\documentclass[reqno]{amsart}
\usepackage{amssymb,setspace}
\usepackage{ifpdf}
\ifpdf
 \usepackage[hyperindex,pagebackref]{hyperref}%
\else
 \expandafter\ifx\csname dvipdfm\endcsname\relax
 \usepackage[hypertex,hyperindex,pagebackref]{hyperref}%
 \else
 \usepackage[dvipdfm,hyperindex,pagebackref]{hyperref}%
 \fi
\fi
\theoremstyle{plain}
\newtheorem{thm}{Theorem}
\numberwithin{equation}{section}
\DeclareMathOperator{\td}{d}
\allowdisplaybreaks[4]

\begin{document}

\title[A double inequality for bounding Toader mean]
{A double inequality for bounding Toader mean by the centroidal mean}

\author[Y.Hua]{Yun Hua}
\address[Hua]{Department of Mathematics, School of Science, Tianjin Polytechnic University, Tianjin City, 300387, China}
\email{\href{mailto: Yun Hua <xxgcxhy@163.com>}{xxgcxhy@163.com}}

\author[F. Qi]{Feng Qi}
\address[Qi]{College of Mathematics, Inner Mongolia University for Nationalities, Tongliao City, Inner Mongolia Autonomous Region, 028043, China}
\email{\href{mailto: F. Qi <qifeng618@gmail.com>}{qifeng618@gmail.com}, \href{mailto: F. Qi <qifeng618@hotmail.com>}{qifeng618@hotmail.com}, \href{mailto: F. Qi <qifeng618@qq.com>}{qifeng618@qq.com}}
\urladdr{\url{http://qifeng618.wordpress.com}}

\begin{abstract}
In the paper, the authors find the best numbers $\alpha$ and $\beta$ such that
$$
\overline{C}\bigl(\alpha a+(1-\alpha)b,\alpha b+(1-\alpha)a\bigr)<T(a,b)
<\overline{C}\bigl(\beta a+(1-\beta)b,\beta b+(1-\beta)a\bigr)
$$
for all $a,b>0$ with $a\ne b$, where $\overline{C}(a,b)={2\bigl(a^2+ab+b^2\bigr)}{3(a+b)}$ and $T(a,b)=\frac{2}{\pi}\int_{0}^{{\pi}/{2}}\sqrt{a^2{\cos^2{\theta}}+b^2{\sin^2{\theta}}}\,\td\theta$
denote respectively the centroidal mean and Toader mean of two positive numbers $a$ and $b$.
\end{abstract}



\thanks{This paper was typeset using \AmS-\LaTeX}

\maketitle

\section{Introduction}

In~\cite{Hua-Qi-toader31}, Toader introduced a mean
\begin{align}\label{eq1.1}
T(a,b)&=\dfrac{2}{\pi}\int_{0}^{\pi/2}\sqrt{a^2\cos^2\theta+b^2\sin^2\theta}\,\td\theta\\
&=
\begin{cases}
\dfrac{2a}\pi\mathcal{E}\Biggl(\sqrt{1-\biggl(\dfrac{b}a\biggr)^2}\,\Biggr),&a>b,\\
\dfrac{2b}\pi\mathcal{E}\Biggl(\sqrt{1-\biggl(\dfrac{a}{b}\biggr)^2}\,\Biggr),&a<b,\\
a,&a=b,
\end{cases}
\end{align}
where
$$
{\mathcal{E}}={\mathcal{E}}(r)=\int_{0}^{\pi/2}\sqrt{1-r^2\sin^2\theta}\,\td\theta
$$
for $r\in[0,1]$ is the complete elliptic integral of the second kind.
\par
In recent years, there have been plenty of literature, such as~\cite{Hua-Qi-toader35, Hua-Qi-toader33, Hua-Qi-toader34, Jiang-toader2.tex, jiang-qi-toader-bound.tex, Jiang-Neuman.tex, jiang-qi-October.tex, background-Jiang-Qi.tex, background-Jiang-Qi-revised.tex, Seiffert-Jiang-Qi.tex, neuman-sandor-mean.tex, Toader-M-Li-Zheng-Rev.tex, Hua-Qi-toader32}, dedicated to bounding Neuman-S\'andor's, Seiffert's, Toader's, and other means which are related to complete elliptic integrals of the second kind.
\par
For $p\in \mathbb{R}$ and $a,b>0$, the centroidal mean $\overline{C}(a,b)$ and the $p$-th power mean $M_{p}(a,b)$ are defined respectively by
\begin{equation}\label{eq1.2}
 \overline{ C}(a,b)=\frac{2\bigl(a^2+ab+b^2\bigr)}{3(a+b)}
\end{equation}
and
\begin{equation}\label{eq1.3}
M_{p}(a,b)=
\begin{cases}
\biggl(\dfrac{a^{p}+a^{p}}{2}\biggr)^{{1}/{p}},&p\ne 0,\\
\sqrt{ab}\,,&p=0.
\end{cases}
\end{equation}
\par
In~\cite{Hua-Qi-toader39}, Vuorinen conjectured that
\begin{equation}\label{eq1.4}
  M_{3/2}(a,b)<T(a,b)
\end{equation}
for all $a,b>0$ with $a\ne b$. This conjecture was verified by Qiu and Shen~\cite{Hua-Qi-toader310} and by Barnard, Pearce, and Richards~\cite{Hua-Qi-toader311}.
In~\cite{Hua-Qi-toader312}, Alzer and Qiu presented that
\begin{align}\label{eq1.5}
  T(a,b)<M_{(\ln2)/\ln(\pi/2)}(a,b)
\end{align}
for all $a,b>0$ with $a\ne b$, which gives a best possible upper bound for Toader mean in terms of the power mean.
\par
Very recently, Chu, Wang, and Ma proved in~\cite{Hua-Qi-toader313} that the double inequality
\begin{equation}\label{eq1.6}
C\bigl(\alpha a+(1-\alpha)b,\alpha b+(1-\alpha)a\bigr)<T(a,b)<C\bigl(\beta a+(1-\beta)b,\beta b+(1-\beta)a\bigr)
\end{equation}
is valid for all $a,b>0$ with $a\ne b$ if and only if $\alpha\le \frac34$ and $\beta\ge \frac12+\frac{\sqrt{4\pi-\pi^2}\,}{2\pi}$, where $C(a,b)=\frac{a^2+b^2}{a+b}$ is the contraharmonic mean.
\par
For positive numbers $a, b>0$ with $a\ne b$, let
\begin{equation}\label{J(x)-dfn-eq}
J(x)=\overline{C}\bigl(xa+(1-x)b,xb+(1-x)a\bigr)
\end{equation}
on $\bigl[\frac{1}{2},1\bigr]$. It is easy to see that $J(x)$ is continuous and strictly increasing on $\bigl[\frac{1}{2},1\bigr]$.
Now it is much natural to ask a question: What are the best constants $\alpha\ge\frac12$ and $\beta\le1$ such that the double inequality
\begin{equation}\label{seiffert-th1.1-constant}
\overline{C}\bigl(\alpha a+(1-\alpha)b,\alpha b+(1-\alpha)a\bigr)<T(a,b)
<\overline{C}\bigl(\beta a+(1-\beta)b,\beta b+(1-\beta)a\bigr)
\end{equation}
holds for $a,b>0$ with $a\ne b$? This problem can be affirmatively answered by the following theorem which is the main result of this paper.

\begin{thm}\label{th1.1}
For positive numbers $a,b>0$ with $a\ne b$, the double inequality~\eqref{seiffert-th1.1-constant} is valid if and only if $\alpha\le \frac{1}{2}\bigl(1+\frac{\sqrt{3}\,}{2}\bigr)$ and $\beta\ge \frac{1}{2}+\frac{1}{2}\sqrt{\frac{12}{\pi}-3}\,$.
\end{thm}

\section{Proof of Theorem~\ref{th1.1}}

For $0<r<1$, denote $r'=\sqrt{1-r^2}\,$. It is known that Legendre's complete elliptic integrals of the first and second kind are defined respectively by
$$
\begin{cases}\displaystyle
{\mathcal{K}}={\mathcal{K}}(r)=\int_{0}^{{\pi}/{2}}\frac1{\sqrt{1-r^2\sin^2\theta}}\td \theta,\\
{\mathcal{K}}'={\mathcal{K}}'(r)={\mathcal{K}}(r'),\\
{\mathcal{K}}(0)=\dfrac\pi2,\\
{\mathcal{K}}(1^-)=\infty
\end{cases}
$$
and
$$
\begin{cases}\displaystyle
{\mathcal{E}}={\mathcal{E}}(r)=\int_{0}^{\pi/2}\sqrt{1-r^2\sin^2\theta}\,\td\theta,\\
{\mathcal{E}}'={\mathcal{E}}'(r)={\mathcal{E}}(r'),\\
{\mathcal{E}}(0)=\dfrac\pi2,\\
{\mathcal{E}}(1^-)=1.
\end{cases}
$$
See~\cite{Hua-Qi-toader36, Hua-Qi-toader37}.
For $0<r<1$, the following formulas were presented in~\cite[Appendix~E, pp.~474\nobreakdash--475]{Hua-Qi-toader38}:
\begin{gather*}
 \frac{\td \mathcal{K}}{\td r}=\frac{\mathcal{E}-(r')^2{\mathcal{K}}}{r(r')^2},\quad
\frac{\td \mathcal{E}}{\td r}=\frac{\mathcal{E}-\mathcal{K}}{r},\quad
\frac{\td ({\mathcal{E}}-(r')^2{\mathcal{K}})}{\td r}=r{\mathcal{K}},\\
\frac{\td ({\mathcal{K}}-{\mathcal{E}})}{\td r}=\frac{r{\mathcal{E}}}{(r')^2},
\quad
{\mathcal{E}}\biggl(\frac{2\sqrt{r}}{1+r}\biggr)=\frac{2\mathcal{E}-(r')^2\mathcal{K}}{1+r}.
\end{gather*}
\par
For simplicity, denote
\begin{equation*}
\lambda=\frac{1}{2}\biggl(1+\frac{\sqrt{3}\,}{2}\biggr)\quad \text{and}\quad \mu=\frac{1}{2}+\frac{1}{2}\sqrt{\frac{12}{\pi}-3}\,.
\end{equation*}
It is clear that, in order to prove the double inequality~\eqref{seiffert-th1.1-constant}, it suffices to show
\begin{equation}\label{eq2.1}
  T(a,b)>\overline{C}\bigl(\lambda a+(1-\lambda)b,\lambda b+(1-\lambda)a\bigr)
\end{equation}
and
\begin{equation}\label{eq2.2}
  T(a,b)<\overline{C}\bigl(\mu a+(1-\mu)b,\mu b+(1-\mu)a\bigr).
\end{equation}
From~\eqref{eq1.1} and~\eqref{eq1.2} we see that both $T(a,b)$ and $\overline{C}(a,b)$ are symmetric and homogenous of degree $1$. Hence, without loss of generality, we assume that $a>b$. Let $t=\frac{b}a\in(0,1)$ and $r=\frac{1-t}{1+t}\in(0,1)$ and let $p\in \bigl(\frac12,1\bigr)$. Then
\begin{multline}\label{eq2.3}
\begin{aligned}
&\quad T(a,b)-\overline{C}\bigl(pa+(1-p)b,pb+(1-p)a\bigr)\\
&=\frac{2a}{\pi}\mathcal{E}\Biggl(\sqrt{1-\biggl(\frac{b}{a}\biggr)^2}\,\Biggr)
\end{aligned}\\
-2a\frac{[p+(1-p)b/a]^2+[p+(1-p)b/a](pb/a+1-p)+(pb/a+1-p)^2}{3(1+b/a)}\\
=\frac{2a}{\pi}\mathcal{E}\bigl(\sqrt{1-t^2}\,\bigr)-2a\frac{[p+(1-p)t]^2 +[p+(1-p)t](pt+1-p)+(pt+1-p)^2}{3(1+t)}\\
\begin{aligned}
&=\frac{2a}{\pi}\frac{2\mathcal{E}-\bigl(1-r^2\bigr)\mathcal{K}}{1+r}-a\frac{(1-2p)^2r^2+3}{3(1+r)}\\
&=\frac{a}{1+r}\biggl\{\frac{2}{\pi}\bigl[2\mathcal{E}-\bigl(1-r^2\bigr)\mathcal{K}\bigr] -\frac{1}{3}(1-2p)^2r^2-1\biggr\}.
\end{aligned}
\end{multline}
Let
\begin{align}\label{eq2.4}
f(r)=\frac{2}{\pi}\bigl[2\mathcal{E}-\bigl(1-r^2\bigr)\mathcal{K}\bigr]-\frac{1}{3}(1-2p)^2r^2-1,
\end{align}
and let $f_1(r)=rf'(r)$ and $f_2(r)=\frac{f_1'(r)}r$. Then, by standard argument, we have
\begin{gather*}
  f(0)=0,\quad f_1(0)=0,\quad f_2(0)=1-\frac{4}{3}(1-2p)^2,\\
  f(1^-)=\frac{4}{\pi}-1-\frac{1}{3}(1-2p)^2,\quad
 f_1(1^-)=\frac{2}{\pi}-\frac{2}{3}(1-2p)^2,\quad   f_2(1^-)=+\infty,\\
f_1(r)=\frac{2}{\pi}\bigl[\mathcal{E}-\bigl(1-r^2\bigr)\mathcal{K}\bigr]-\frac{2}{3}(1-2p)^2r^2,\quad
   f_2(r)=\frac{2}{\pi}\mathcal{K}-\frac{4}{3}(1-2p)^2,
\end{gather*}
\par
When $p=\lambda=\frac{1}{2}\bigl(1+\frac{\sqrt{3}\,}{2}\bigr)$, it follows that $f_2(0)=0$. An easy argument leads to $f(r)>0$ for $r\in(0,1)$. Together with this, the inequality~\eqref{eq2.1} follows from~\eqref{eq2.3} and~\eqref{eq2.4}.
\par
When $p=\mu=\frac{1}{2}+\frac{1}{2}\sqrt{\frac{12}{\pi}-3}\,$, it is simple to derive that
$$
f(1^-)=0,\quad f_1(1^-)=\frac{2(\pi-3)}{\pi}>0,\quad f_2(0)=\frac{5\pi-16}{\pi}<0.
$$
Consequently, considering the monotonicity of $f_2(r)$, it is deduced that there exists $r_0\in(0,1)$ such that $f_2(r)<0$ on $(0,r_0)$ and $f_2(r)>0$ on $(r_0,1)$. Hence, the function $f_1(r)$ is strictly decreasing on $(0,r_0)$ and strictly increasing on $(r_0,1)$.
Similarly, there exists $r_1\in(0,1)$ such that $f_1(r)<0$ on $(0,r_1)$ and $f_1(r)>0$ on $(r_1,1)$. Thus, the function $f(r)$ is strictly decreasing on $(0,r_1)$ and strictly increasing on $(r_1,1)$.
As a result, the inequality~\eqref{eq2.2} follows.
\par
If $p>\lambda$, then $f_2(r)<0$. From the continuity of $f(r)$, $f_1(r)$, and $f_2(r)$, it follows that there exists $\delta_1=\delta_1(p)>0$ such that $f(r)<0$ on $(0,\delta_1)$. Combining this with~\eqref{eq2.3} and~\eqref{eq2.4} yields that $T(a,b)<\overline{C}\bigl(p a+(1-p)b,p b+(1-p)a\bigr)$ for $\frac{b}a\in\bigl(\frac{1-\delta_1}{1+\delta_1},1\bigr)$.
If $p<\mu$, then $f(1^-)>0$. Hence, there exists $\delta_2=\delta_2(p)\in(0,1)$ such that $f(r)>0$ on $(1-\delta_2,1)$. Combining this with~\eqref{eq2.3} and~\eqref{eq2.4} reveals that $T(a,b)>\overline{C}\bigl(p a+(1-p)b,p b+(1-p)a\bigr)$ for $\frac{b}a\in(0,\delta_2/(2-\delta_2))$. These imply that the constants $\lambda$ and $\mu$ are the best possible.
The proof of Theorem~\ref{th1.1} is complete.

\subsection*{Acknowledgements}
The authors thank the anonymous referee for his/her careful reading and helpful corrections to the original version of this paper.
\par
The first author was partially supported by the Project of Shandong Province Higher Educational Science and Technology Program under Grant No. J11LA57, China.

\end{document}